\documentclass[a4page]{amsart}

\newcommand{\euler}{\mathrm{e}}
\newcommand{\uc}{\mathrm{uc}}

\newcommand{\RR}{\mathbb{R}}

\newcommand{\NN}{\mathbb{N}}
\newcommand{\PP}{\mathbb{P}}
\newcommand{\EE}{\mathbb{E}}
\newcommand{\ZZ}{\mathbb{Z}}
\newcommand{\equiset}{S_{\delta , Z}}
\newcommand{\ran}{\operatorname{Ran}}

\newtheorem{theorem}{Theorem}

\begin{document}

\title{Uncertainty relations and applications in spectral and control theory}
\author{Ivan Veseli\'c}

\begin{abstract}

This is an extended abstract for the Oberwolfach Workshop 2101b
\emph{Geometry, Dynamics and Spectrum of Operators on Discrete Spaces}
organized by
David Damanik, Matthias Keller, Tatiana Smirnova-Nagnibeda, and Felix Pogorzelski.
It took place in online format from 3rd to 9th January 2021.
\end{abstract}

\maketitle

Unique continuation is a prominent phenomenon encountered in the study of several
classes of functions, e.~g.~subsets of holomorphic functions or spaces of solutions of partial
differential equations.
The study of this phenomenon, in the multidimensional setting, goes back at least to
Carleman and M\"uller. It is impossible to give here an adequate overview
of the gradual development of understanding of this phenomenon.
Let us just mention that it has crucial consequences in a plethora of applied problems:
Absence of eigenvalues embedded in the continuous spectrum of Schr\"odinger operators,
size and dimensionality of nodal sets of Laplace-Beltrami operators on compact manifolds,
inverse problems for partial differential equations,
like the Calderon problem or the control problem for evolution equations.

More recently,
the unique continuation principle has been successfully applied
to the study of mathematical models in solid state and condensed matter physics.
A particular feature of these models is the presence of a microscopic and
a macroscopic length scale, where the latter is many orders of magnitude larger than the former.
For periodic and ergodic Schr\"odinger operator,
the ratio between the macroscopic and microscopic scale is determined by the spacing between atoms in crystals,
which is typically a few \r{A}ngstr\"oms.

Consequently, the unique continuation principles have to take into account this
geometric structure to be applicable in the mentioned physical context.
This means that on the technical level we have to consider the unique continuation
problem on large domains --- `large' compared to some reference scale,
e.~g.~the atomic scale in the case of ergodic Schr\"odinger operators.
It is this reference scale at which changes of the coefficient functions
(e.~g.~the electric potential) are observed.

As in other situations of statistical and condensed matter physics,
it is natural to approximate large bounded domains by unbounded ones,
and hence include the latter in our analysis.
Now we have set the stage to formulate unique continuation
estimates on unbounded and (large) bounded domains
and spell out thereafter applications in three different questions of mathematical physics.

\section{Quantitative scale free unique continuation principles}

Let $d \in \NN$, $G > 0$, $\delta > 0$  and $\Gamma = 
           \times_{i =1}^d (\alpha_i , \beta_i) \subset \RR^d$ with
$\alpha_i ,\beta_i \in \RR\cup \{\pm \infty\}$. Assume that  $\Lambda_G:=(-G/2 , G/2)^d \subset \Gamma$. We say that a sequence $Z = (z_j)_{j \in (G\ZZ)^d } \subset \RR^d$ is \emph{$(G,\delta)$-equidistributed}, if
 \[
  \forall j \in (G\ZZ)^d  \colon \quad  B(z_j , \delta) \subset (-G/2 , G/2)^d + j.
\]
For a $(G,\delta)$-equidistributed sequence $Z$ define
\[
\equiset = \bigcup_{j \in (G\ZZ)^d } B(z_j , \delta) \cap \Gamma.
\]
For a real $V \in L^\infty(\Gamma)$ define $ H = -\Delta + V$ on $L^2 (\Gamma)$ with Dirichlet or Neumann boundary conditions.
\smallskip
\begin{theorem}[\cite{NTTV-JST}] \label{thm:sfucp}
There is $N>0$ depending only on $d$, such that for all
$G>0$, all $\Lambda_G\subset\Gamma \subset\RR^d$ as above, all $\delta \in (0,G/2)$, all $(G,\delta)$-equidistributed sequences $Z$, all  $V \in L^\infty (\Gamma)$, all $E \in \RR$, and all $\psi \in \ran \mathbf{1}_{(E , \infty)}(H)$ we have
\begin{equation*}
\lVert \psi \rVert_{L^2 (\equiset)}^2
\geq C_\uc(V,E) \lVert \psi \rVert_{L^2 (\Gamma)}^2 ,
\quad\text{where}\quad
t_+:=\max\{0,t\} \text{ for } t \in \RR
 \end{equation*}
and $\displaystyle C_\uc(V,E) = \sup_{\lambda \in \RR}
 \left(\frac{\delta}{G}\right)^{N \bigl(1 + G^{4/3}\lVert V-\lambda \rVert_\infty^{2/3} + G\sqrt{(E-\lambda)_+} \bigr)}
$.
\end{theorem}
The estimate is called \emph{scale free} since it is independent of the domain $\Gamma$.

\section{Lifting estimates of edges of essential spectrum by potentials}
The min max principle for hermitean matrices shows and quantifies the lifting of eigenvalues
under the influence of a positive definite perturbation.
In the case of Schr\"odinger operators, eigenvalues may have accumulation points and/or be located inside gaps of the essential spectrum.
The question is whether analogous shifting estimates hold in these cases.
Furthermore:
Do similar lifting estimates hold for the edges of essential spectrum;
and: Are positive semi-definite potentials sufficient to produce such lifting?
We have a positive answer to these questions formulated in the following

\begin{theorem}[\cite{NTTV-JST}] \label{thm:lifting-essential}
Let $N$, $G$, $\Gamma $, $\delta$, $Z$, $V$, $H$, and $C_\uc(V,E)$ be as in
Theorem \ref{thm:sfucp}. Let $W \in L^\infty (\Gamma)$ be real-valued with $ W \geq \vartheta \mathbf{1}_{\equiset}$ for some $\vartheta > 0$.

 Let $a,b \in \sigma_{\mathrm{ess}}(H)$, and $a< b$ such that $(a,b) \cap \sigma_{\mathrm{ess}}(H) = \emptyset$. We set $t_0 = (b-a)/\lVert W \rVert_\infty$,  and $t_+= \max\{0, t\}$, $t_-= \max\{0, -t\}$.
Then the functions $f_\pm : (-t_0 , t_0) \to \RR$
\begin{align*}
 f_-(t) &= \sup \left( \sigma_{\mathrm{ess}}(H + t W) \cap (-\infty, b - t_- \lVert W \rVert_\infty) \right), \\
 f_+ (t) &= \inf \left( \sigma_{\mathrm{ess}}(H + t W) \cap (a + t_+ \lVert W \rVert_\infty  , \infty) \right) ,
\end{align*}
satisfy for all $t \in (-t_0,t_0), \varepsilon >0$ such that $t+\varepsilon \in (-t_0,t_0)$
the two-sided Lipschitz bound
\begin{align*}
\varepsilon \, \vartheta \, \sup_{0\leq t\leq1}C_\uc(V+tW,b + \|W\|_\infty) &\leq f_\pm(t+\varepsilon)-f_\pm(t)
\leq \varepsilon \, \Vert W\Vert
\end{align*}
\end{theorem}
Analogous bounds hold for discrete eigenvalues in all gaps of the essential spectrum (and below it).

\section{Anderson localization for random Schr\"odinger operators}
It is well known from previous work that unique continuation estimates
are a powerful tool for deriving \emph{Wegner estimates}, which in turn play a crucial role in proving localization. We present here a result for a model with non-linear parameter dependence.

For $0 \leq \omega_{-} < \omega_{+} < \tfrac{1}{4}$ set
$ \Omega = \times_{j \in {\ZZ^d}} \RR$,
$ \PP = \bigotimes_{j \in {\ZZ^d}} \mu$
where
$\mu$ is a probability measure with $\mathrm{supp}\ \mu \subset [\omega_{-}, \omega_{+}]$ and a bounded density $\nu$.
Hence, $\pi_j(\omega) \mapsto \omega_j,j \in {\ZZ^d}$, are continuous iid random variables.
The standard random breather model is defined as
\begin{equation}\label{eq:standardRBP}
H_\omega = -\Delta + V_\omega^{\mathrm{br}}(x) \quad \text{with} \quad V_\omega^{\mathrm{br}}(x) = \sum_{j \in \ZZ^d} \chi_{B_{\omega_j}}(x-j),
\end{equation}
and its restriction to the box $\Lambda_L$ (with Dirichlet or Neumann boundary conditions) is denoted by $H_{\omega,L}$.
\begin{theorem}[Wegner estimate for the standard random breather model
\cite{TaeuferV-15,NTTV-APDE}]\label{thm:Wegner}
Fix $E_0 \in \RR$
and set $ \varepsilon_{\max} = (1/4)\cdot  8^{-N(2+{\lvert E_0+1 \rvert}^{1/2})}$,
where $N$ is the constant from Theorem~\ref{thm:sfucp}.
Then there is  $C=C(d,E_0) \in (0,\infty)$ such that for all $\varepsilon \in (0,\varepsilon_{\max}]$ and $E \geq 0$ with
$[E-\varepsilon, E+\varepsilon] \subset (- \infty , E_0]$,
we have
\begin{equation*}
\EE \left[ \mathrm{Tr} \left[ \chi_{[E- \varepsilon, E + \varepsilon]}(H_{\omega,L}) \right] \right]
\leq
C
\lVert \nu \rVert_\infty
\varepsilon^{[N(2+{\lvert E_0 + 1 \rvert}^{1/2})]^{-1}}
\left\lvert\ln \varepsilon \right\rvert^d L^d.
\end{equation*}
\end{theorem}

Of similar importance for the proof of localization ist the so-called
\emph{initial scale estimate}. {Seelmann} and {T\"aufer} succeeded in proving one (and consequently localization) for energies near spectral band edges of randomly perturbed periodic potentials.
It is remarkable that they do not need to assume anything about the extrema for the associated Floquet eigenvalues, nor on the vanishing order of $\nu$ near the extrema of $supp \nu$
(in contrast to \cite{BarbarouxCH-97b,KirschSS-98a,Klopp-99,Veselic-02b}).

\begin{theorem}[Localization at band edges \cite{ST}]\label{thm:ST-loc}
Let $V_{\mathrm{per}}\colon \RR^d \to \RR$ be bounded and periodic, $H_{\mathrm{per}}=-\Delta+V_{\mathrm{per}}$,
$V_\omega^{\mathrm{an}}= \sum_{j\in\ZZ^d} \omega_j u(\cdot-j)$ with iid bounded random variables $\omega_j$ with bounded density $\nu$ and
$L^\infty_c(\RR^d)\ni u\geq c \, \chi_{B_\delta}$, for some $c,\delta>0$, and $H_\omega=H_{\mathrm{per}}+V_\omega^{\mathrm{an}}$.
Assume that $a< b$, $(a,b)\subset \rho(H_\omega)$, and $b \in \sigma(H_\omega)$. Then there is an $\varepsilon>0$ such that in $[b,b+\varepsilon]$ the operator $H_\omega$ exhibits Anderson localization.
\end{theorem}

\section{Null control and observability estimates for the heat equation}
The following observability (and hence the associated control cost) estimates for the generalized heat equation (on $\Gamma$ as above) have been derived from Theorem \ref{thm:sfucp}.

\begin{theorem}[\cite{NTTV-ESAIM}]\label{thm:control}
Let $N$, $G$, $\Gamma $, $\delta$, $Z$, $V$, and  $H$ be as in
Theorem \ref{thm:sfucp}.
Then for all $\phi \in L^2 (\Gamma)$, and all $T > 0$ we have
\begin{align*}
  \lVert \euler^{- HT} \rVert_{L^2 (\Gamma)}^2
  &\leq
  C_\mathrm{obs}(\delta,G,\|V\|_\infty,T)^2
  \int_0^T \lVert \euler^{- Ht} \phi \rVert_{L^2 (\equiset)}^2 \mathrm{d} t ,
\end{align*}
where for some $C_1, C_2,C_3>0$ depending only on the dimension we have
\begin{align*}
 C_\mathrm{obs}
 &=  \left( \frac{\delta}{G} \right)^{-C_2 (1 + G^{4/3} \lVert V \rVert_\infty^{2/3})} \frac{C_1}{T} \exp \left( \frac{C_3 G^2 \ln^2 (\delta / G)}{T} \right) \ \text{ if $\kappa:=\inf \sigma(H)  \geq 0$}, \\[2ex]
 C_\mathrm{obs}
&=  \left( \frac{\delta}{G} \right)^{-C_2 (1 + G^{4/3} \lVert V - \kappa \rVert_\infty^{2/3})} \inf_{t \in [0,T)}
      \frac{C_1}{T-t}
      \exp \left( \frac{C_3 G^2 \ln^2 (\delta / G)}{T-t} - 2 \kappa t \right)
\end{align*}
if $\kappa>0$.
\end{theorem}

%

\def\polhk#1{\setbox0=\hbox{#1}{\ooalign{\hidewidth
  \lower1.5ex\hbox{`}\hidewidth\crcr\unhbox0}}}

\end{document}